\documentclass[12pt]{article}

\usepackage{amsmath,enumerate,amssymb,amsthm}
\newtheorem{Theorem}{Theorem}

\theoremstyle{definition}
\newtheorem{defn}{Definition}

\begin{document}

\title{Divisor Functions and the Number of Sum Systems}
 \author{{Matthew C. Lettington} and {Karl Michael Schmidt}}
\maketitle
\begin{abstract}

\noindent
Divisor functions have attracted the attention of number theorists from Dirichlet to the present day.
Here we consider associated divisor functions
$c_j^{(r)}(n)$
which for non-negative integers
$j, r$
count the number of ways of representing
$n$
as an ordered product of
$j+r$
factors, of which the first
$j$
must be non-trivial, and their natural extension to negative integers
$r.$
We give recurrence properties and explicit formulae for these novel arithmetic functions.
Specifically, the functions
$c_j^{(-j)}(n)$
count, up to a sign, the number of ordered factorisations of
$n$
into
$j$
square-free non-trivial factors. These functions are related to a modified version of the M\"obius function and
turn out to play a central role in counting the number of sum systems of given dimensions.
\par
Sum systems are finite collections of finite sets of non-negative integers, of prescribed cardinalities, such that their set sum generates consecutive integers without repetitions. Using a recently
established bijection between sum systems and joint ordered factorisations of their component set cardinalities,
we prove a formula expressing the number of different sum systems in terms of associated divisor functions.
\end{abstract}


\def\rc#1{\frac{1}{#1}}
\def\qed{\square}
\def\cases#1{\left\{\begin{matrix} #1 \end{matrix}\right.}
\def\pmatrix#1{\left(\begin{matrix} #1 \end{matrix}\right)}
\def\nonpmatrix#1{\begin{matrix} #1 \end{matrix}}


\section{Introduction} \label{s1}
\par\medskip\noindent
The recent work \cite{rHLS} gives a construction formula for all sum systems,
which are finite collections of finite sets of non-negative integers, of
prescribed cardinalities, such that the process of taking one element from
each component set and adding them up generates each number in an arithmetic
progression exactly once.
Thus a sum system is of the form
$A_1, \dots, A_m \subset {\mathbb N}_0,$
\begin{align}
 \sum_{k=1}^m A_k &= \left\{0, 1, \dots, \prod_{k=1}^m |A_k| - 1 \right\}.
\nonumber\end{align}
Such systems play a role in the construction of principal reversible cuboids
and (in the simple case of only two component sets) of other matrices with
integer entries and
specific symmetry properties, e.g.\ most-perfect pandiagonal squares, cf.\
\cite{rOB}.
The construction of sum systems of given cardinalities
$a_1, \dots, a_m$
for the
$m$
component sets is based on a joint ordered factorisation of these
cardinalities, defined as follows (cf. \cite{rHLS}, Definition 6.6).
\par\medskip\noindent
\begin{defn}
Let
$m \in {\mathbb N}$
and
$a \in {\mathbb N}^m.$
Then we call
\begin{align}
 &((j_1, f_1), (j_2, f_2), \dots, (j_L, f_L)) \in (\{1, \dots, m\} \times ({\mathbb N}+1))^L,
\nonumber\end{align}
where
$L \in {\mathbb N},$
a
{\it joint ordered factorisation\/}
of
$a = (a_1, \dots, a_m)$
if
\begin{align}
 \prod_{j_l = j} f_l = a_j &\qquad (j \in \{1, \dots, m\})
\nonumber\end{align}
and
$j_l \neq j_{l-1}$
$(l \in \{2, \dots, L\}).$
\end{defn}
In other words, a joint ordered factorisation of
an $n$-tuple of natural numbers
$a_1, \dots, a_m$
arises from writing each of these numbers as a product of non-trivial factors,
i.e.\ factors~$\ge~2$, and then arranging all factors in a linear chain such
that no two adjacent factors arise from the factorisation of the same number.
Given a joint ordered factorisation, the sets
\begin{align}
 A_k &= \sum_{j_l = k} \left(\prod_{s=1}^{l-1} f_s \right) \{0, 1, \dots, f_l-1\}
 \qquad (j \in \{1, \dots, m\})
\nonumber\end{align}
form a sum system, and conversely any sum system arises from some joint ordered
factorisation of its dimensions in this way (cf.\ \cite{rHLS}, Theorem 6.7).
This establishes a bijection between sum systems and joint ordered factorisations.
\par
As an illustrative example, consider the case
$m = 5,$
$(a_1, a_2, a_3, a_4, a_5) = (4, 6, 8, 12, 20).$
An example of a joint ordered factorisation of this quintuple of dimensions is
\begin{align}
 &((1, 2), (5, 2), (2, 2), (5, 5), (3, 4), (5, 2), (4, 4), (3, 2), (4, 3), (1, 2), (2, 3)),
\nonumber\end{align}
yielding the corresponding sum system
\begin{align}
 A_1=&\{0,1,7680,7681\}
\nonumber\\
 A_2=&\{0, 4, 15360, 15364, 30720, 30724\}
\nonumber\\
 A_3=&\{0, 40, 80, 120, 3840, 3880, 3920, 3960\}
\nonumber\\
 A_4=&\{0, 320, 640, 960, 1280, 1600, 1920, 2240, 2560, 2880, 3200, 3520\}
\nonumber\\
 A_5=&\{0, 2, 8, 10, 16, 18, 24, 26, 32, 34, 160, 162, 168, 170, 176, 178, 184, 186, 192, 194\}.
\nonumber\end{align}
Thus, inscribing the above numbers on the faces of the five platonic polyhedra
and adding the numbers obtained in each roll of these five dice, we obtain a
random number generator for the integers
$0, 1, \dots, 46079\ (= a_1 a_2 a_3 a_4 a_5-1)$
with uniform probability distribution.
\par
In the present paper, we answer the question of how many different joint
ordered factorisations of a given $m$-tuple of positive integers there are.
Our main result is Theorem~\ref{tcountjof}, which expresses the number of joint ordered
factorisations in terms of values of certain associated divisor functions at
$a_1, \dots, a_m$
(see Eq.\
$(\ref{eNa})$
).
These functions turn out to be closely linked to a modified version of the
number theoretic M\"obius function (cf. \cite{rGEA} p.\ 77). In Section 2, we study some of their properties
before tackling the counting problem in Sections 3 and 4.\vskip 30pt

\section{Non-trivial and Associated Divisor Functions}\label{sprelim}

\noindent
Divisor functions have been studied by many eminent number theorists, from Dirichlet to the present day (e.g.\ \cite{rRama}, \cite{rIvic}, \cite{rKRRR}).
The non-trivial and associated divisor functions defined in \cite{rHHLS} can be
conveniently described in the framework of the commutative Dirichlet convolution algebra
of arithmetic functions.
The convolution of arithmetic functions
$f_1, f_2, \dots, f_j$
is given by
\begin{align}
 (f_1 * f_2 * \cdots * f_j)(n) &= \sum_{n_1 n_2 \cdots n_j = n} f_1(n_1) f_2(n_2) \cdots f_j(n_j),
\label{econvo}\end{align}
summing over all ordered factorisations of
$n \in {\mathbb N}$
into
$j$
factors.
We denote the $j$th convolution power as follows,
$f^{*j} := f * f * \cdots * f,$
where the right-hand side has
$j$
repetitions of
$f;$
by the usual convention,
$f^{*0} = e.$
The function
$e(n) = \delta_{n, 1}$
$(n \in {\mathbb N})$
is the neutral element of the Dirichlet convolution product, and the
convolution inverse of the constant function 1 is the well-known M\"obius
function
$\mu$.
\par
In analogy to the standard $j$th divisor function
$d_j = 1^{*j}$
(cf. \cite{rSchSp} p.\ 9),
which counts the ordered factorisations of its argument into
$j$
positive integer factors, we define the
{\it $j$th non-trivial divisor function\/}
$c_j = (1-e)^{*j},$
which counts the ordered factorisations of its argument into
$j$
non-trivial integer factors, i.e.\ into factors
$> 1.$
\par
Furthermore, for non-negative integer
$r,$
the
{\it associated $(j, r)$-divisor function\/}
is defined as
$c_j^{(r)} = (1-e)^{*j} * 1^{*r}.$
In view of Eq.\
$(\ref{econvo}),$
it counts the ordered factorisations of its argument into
$j+r$
factors, of which the first
$j$
must be non-trivial.
\par
As the constant function 1 has a convolution inverse, the latter definition extends
naturally to negative upper indices, giving the associated $(j, -r)$-divisor
function
$c_j^{(-r)} = (1-e)^{*j} * \mu^{*r}.$
(Note that
$1-e$
does not have a convolution inverse, so there is no analogous extension to
negative lower indices.)
The functions
$c_0^{(-r)} = \mu^{*r}$
were studied by Popovici \cite{rPopo}.
In the associated $(j, -r)$-divisor functions, the
{\it modified M\"obius function\/}
\begin{align}
 (\mu - e)(n) &= \cases{ (-1)^{\Omega(n)} & \hbox{\rm if $n$ is square-free} \cr 0 & \hbox{\rm otherwise (including the case $n = 1$)}}
 \qquad (n \in {\mathbb N}),
\nonumber\end{align}
where
$\Omega(n)$
is the number of prime factors of
$n,$
appears naturally. Indeed, if
$j \ge r,$
then
\begin{align}
 c_j^{(-r)} &= (1-e)^{*j-r} * ((1-e) * \mu)^{*r} = (-1)^r (1-e)^{*j-r} * (\mu-e)^{*r};
\nonumber\end{align}
if
$j < r,$
then
\begin{align}
 c_j^{(-r)} &= ((1-e) * \mu)^{*j} * \mu^{*r-j} = (-1)^j (\mu-e)^{*j} * \mu^{*r-j}.
\nonumber\end{align}
Note that
$c_j^{(r)}(n)$
involves factorisation of
$n$
into
$j+r$
factors if
$r \ge 0,$
into
$\max\{j, -r\}$
factors if
$r < 0,$
of which at least
$j$
must be non-trivial, so
$c_j^{(r)}(n) = 0$
if
$j > \Omega(n).$
(Also, if $r < 0$, then at least $-r$ factors must be square-free.)

\par
The special case
$j = -r,$
\begin{align}
 c_j^{(-j)}(n) &= (-1)^j \sum_{n_1 n_2 \cdots n_j = n} (\mu-e)(n_1)\, (\mu-e)(n_2) \cdots (\mu-e)(n_j) \qquad (n \in {\mathbb N}),
\label{especial}\end{align}
turns out to be of particular importance (cf.\ Theorem \ref{tcountjof} below).
The value of
$c_j^{(-j)}(n)$
can be interpreted as
$(-1)^{\Omega(n) + j}$
times the number of ordered factorisations of
$n$
into
$j$
non-trivial, square-free factors.
\par
The following statement aids the calculation of the associated divisor functions
either via a recurrent scheme similar to Pascal's triangle, or
directly in terms of the prime factorisation of their argument.
\begin{Theorem}\label{tadifu}
Let
$j \in {\mathbb N}_0,$
$r \in {\mathbb Z}.$
Then
\begin{description}
\item{(a)}
\begin{align}
 c_j^{(r+1)} &= c_{j+1}^{(r)} + c_j^{(r)};
\label{erecur}\end{align}
\item{(b)}
if
$n = p_1^{a_1} p_2^{a_2} \cdots p_\nu^{a_\nu}$
with distinct primes
$p_1, p_2, \dots, p_\nu,$
then
\begin{align}
 c_j^{(r)}(n) &= \sum_{k=0}^j (-1)^k \pmatrix{j \cr k} \prod_{l=1}^\nu \pmatrix{a_l + r + j - k - 1 \cr a_l}.
\label{efactor}\end{align}
\end{description}
\end{Theorem}
\par\medskip\noindent
{\it Proof.\/}
Eq.\ (\ref{erecur}) follows immediately from the observation that
$c_{j+1}^{(r)} = (1-e)^{*j+1} * 1^{*r} = (1-e)^{*j} * 1 * 1^{*r} - (1-e)^{*j} * e* 1^{*r}.$
For part (b), the binomial theorem gives
\begin{align}
 c_j^{(r)} &= \sum_{k=0}^j (-1)^k \pmatrix{j \cr k} e^{*k} * 1^{*j-k+r},
\label{erel3}\end{align}
and Eq.\
$(\ref{efactor})$
follows from the identity (cf. \cite{rHHLS} Lemma 1)
\begin{align}
 1^{*j}(n) &= \prod_{k=1}^\nu \pmatrix{a_k + j - 1 \cr a_k},
\label{ediv}\end{align}
which holds for all integers
$j.$
As
$1$
and
$1^{*-1} = \mu$
are multiplicative arithmetic functions and the Dirichlet convolution of
multiplicative functions is multiplicative, it is sufficient to verify
$(\ref{ediv})$
for a single prime power.
For positive
$j,$
$1^{*j}(p^a)$
is, by Eq.\
$(\ref{econvo}),$
equal to the number of $j$-part partitions of
$a,$
i.e.\ to
$\pmatrix{a + j - 1 \cr a}.$
Furthermore, again by Eq.\
$(\ref{econvo}),$
$1^{*-j}(p^a) = \mu^{*j}(p^a)$
is equal to
$(-1)^a$
times the number of ways of writing
$a$
as an ordered sum of
$j$
terms, each either
$0$
or
$1,$
i.e.\ to
$(-1)^a \pmatrix{j \cr a} = \pmatrix{-j + a - 1 \cr a}.$
\phantom{.}\hfill$\qed$\par\noindent
\par\medskip\noindent
In the specific case of a power of a square-free number $n$, the product in Eq.\
$(\ref{efactor})$
turns into a power; then, using the last identity in the above proof, we can
derive the formula
\begin{align}
 c_j^{(r-j+1-a)}(n^a) &= (-1)^{a \Omega(n) + j} c_j^{(-r)}(n^a).
\nonumber\end{align}

We note the following relationships between the associated divisor functions
and the standard divisor functions $d_j$ (and their inverses with respect to
Dirichlet convolution).
\begin{Theorem}\label{trelations}
Let
$j \in {\mathbb N}_0,$
$r \in {\mathbb Z}.$
Then
\begin{align}
 d_r &= \sum_{k=0}^\infty \pmatrix{k+r-1 \cr k} c_k^{(-k)}.
\label{erel1}\end{align}
More generally, for any
$u \in {\mathbb N}_0$
and
$v \in {\mathbb Z},$
\begin{align}
 c_{j+u}^{(r+v)} &= \sum_{k=j}^\infty \pmatrix{k+r-1 \cr k-j} c_{u+k}^{(v-k)}.
\label{erel2}\end{align}
\end{Theorem}
\par\medskip\noindent
{\it Proof.\/}
The identity
$(\ref{erel1})$
follows from the inverse binomial formula,
\begin{align}
 1^{*r} &= (e + \mu - e)^{*-r}
 = \sum_{k=0}^\infty (-1)^k \pmatrix{r+k-1 \cr k} (\mu-e)^{*k}.
\nonumber\end{align}
The series is pointwise convergent because
$\mu - e$
is pointwise nilpotent in the convolution algebra.
Hence induction on
$j$
gives
\begin{align}
 (1-e)^{*j} * 1^{*k} &= \sum_{k=j}^\infty \pmatrix{k+r-1 \cr k-j} (e - \mu)^{*k},
\nonumber\end{align}
and Eq.\
$(\ref{erel2})$
follows by convolution with
$(1-e)^{*u}*1^{*v}$
on both sides.
\hfill$\qed$\par\noindent
\par\medskip\noindent
{\it Remark.\/}
Curiously, the binomial coefficient appearing in Eqs
$(\ref{erel1})$
and
$(\ref{erel2})$
can be expressed as the associated divisor function of a $k$th prime power,
$\displaystyle c_j^{(r)}(p^k) = \pmatrix{k+r-1 \cr k-j}$
(cf. \cite{rHHLS} Lemma 5, \cite{rSchSp} p. 62 for $j=0$), giving an alternate form as the sum over products of the form $c_j^{(r)}(p^k)c_k^{(-k)}(n)$.
For
$j = 0,$
this is equal to the number of weak compositions of
$k$
into
$r$
parts (cf.\ \cite{rStan} p.\ 15).
We also note that Eq.\
$(\ref{erel3})$
provides a converse to Eq.\
$(\ref{erel1}).$
\par
Taking
$r = 2$
in Eq.\
$(\ref{erel1})$
gives an expression for the standard divisor function (number of divisors),
$d_2 = \sum \limits_{k=0}^\infty (k+1) c_k^{(-k)}.$
Taking
$r = 1$
yields the identity
$1 = \sum \limits_{k=0}^\infty c_k^{(-k)}.$
We note that the sum
$\sum \limits_{k=0}^\infty |c_k^{(-k)}|,$
which gives the number of ordered factorisations into (any number of) square-free, non-trivial factors, generates all odd integers; indeed,
$\sum \limits_{k=0}^\infty |c_k^{(-k)}(p_1 p_2^m)| = 2 m + 1$
$\ (m \in {\mathbb N}_0).$\vskip 30pt

\section{An Auxiliary Counting Problem}\label{s2}
\par\medskip\noindent
We now turn to the question of counting the number of joint ordered factorisations of a given $m$-tuple. In the present section,
we first consider the following combinatorial problem.
\it
Given a number of coloured (but otherwise identical) blocks, with any number of blocks to each of several colours,
in how many ways can all blocks be arranged in a linear sequence such that no
two adjacent blocks have the same colour?
\rm
Note that the answer may very well be 0; indeed, if there are 2 more blocks of
one colour than of all the other colours taken together, then there is no
possible arrangement.
\par
In the following, we make extensive use of the standard multi-index notation summarised in
the Appendix below.
We denote by
$e_n$
the number of different ways
$|n|$
objects, of which there are
$n_j$
of type
$j,$
$j \in \{1, \dots, m\},$
and which are otherwise indistinguishable,
can be linearly arranged such that no neighbouring objects have the same type.
Then, for any
$n \in {\mathbb N}^m,$
the identity
\begin{align}
 \pmatrix{|n| \cr n} &= \sum_{1_m \le k \le n} \pmatrix{n - 1_m \cr k - 1_m} e_k
\label{ecollapse}\end{align}
holds.
Indeed,
$\pmatrix{|n| \cr n}$
is the number of linear arrangements of all objects ignoring the non-adjacency
condition.
Given any such arrangement, consider the associated collapsed arrangement where
any group of contiguous objects of the same type is replaced with a single such
object, resulting in an arrangement of size
$k \le n$
satisfying the adjacency condition.
There are
$e_k$
different collapsed arrangements of size
$k,$
and
$\pmatrix{n - 1_m \cr k - 1_m}$
different arrangements giving rise to each collapsed arrangement.
\begin{Theorem}\label{tcolours}
Let
$n \in {\mathbb N}^m,$
$m \in {\mathbb N}.$
Then
\begin{align}
 e_n &= \sum_{0_m \le k \le n - 1_m} (-1)^{|k|} \pmatrix{n - 1_m \cr k} \pmatrix{|n-k| \cr n - k}.
\label{eeformula}\end{align}
\end{Theorem}
\par\medskip\noindent
{\it Proof.\/}
The power series for the generating function
\begin{align}
 &\sum_{k \in {\mathbb N}_0^m} e_{1_m+k}\,\frac{x^k}{k!}
\nonumber\end{align}
can be shown to be convergent for all
$x \in {\mathbb R}^m$
by comparison with exponential series using the a priori estimate
\begin{align}
 e_{1_m+k} &\le \pmatrix{|k+1_m| \cr k + 1_m} \le m^{|k|+m},
\nonumber\end{align}
where we used Eq.\
$(\ref{ecollapse})$
in the first and the multinomial theorem in the second inequality.
Using the exponential series
\begin{align}
 \exp \sum_{j=1}^m x_j &= \sum_{k \in {\mathbb N}_0^m} \frac{x^l}{k!}
 \qquad (x \in {\mathbb R}^m)
\nonumber\end{align}
and the identity
$(\ref{ecollapse})$
between two applications of the multivariate Cauchy product formula, we find
\begin{align}
 \sum_{k \in {\mathbb N}_0} e_{1_m+k}\,\frac{x^k}{k!} &=
 \left(\sum_{k \in {\mathbb N}_0} e_{1_m+k}\,\frac{x^k}{k!}\right) \left(\sum_{l \in {\mathbb N}_0^m} \frac{x^l}{l!}\right) \left(\sum_{l \in {\mathbb N}_0^m} (-1)^{|l|} \frac{x^l}{l!}\right)
\nonumber\\
 &= \left(\sum_{n \in {\mathbb N}_0} \left(\sum_{0_m \le k \le n} \pmatrix{n \cr k} e_{1_m+k}\right) \frac{x^n}{n!}\right) \left(\sum_{l \in {\mathbb N}_0^m} (-1)^{|l|} \frac{x^l}{l!}\right)
\nonumber\\
 &= \left(\sum_{n \in {\mathbb N}_0} \pmatrix{|n+1_m| \cr n+1_m} \frac{x^n}{n!}\right) \left(\sum_{l \in {\mathbb N}_0^m} (-1)^{|l|} \frac{x^l}{l!}\right)
\nonumber\\
 &= \sum_{k \in {\mathbb N}_0^m} \left(\sum_{0_m \le l \le k} (-1)^{|l|} \pmatrix{k \cr l} \pmatrix{|k+1_m-l| \cr k+1_m-l} \right) \frac{x^k}{k!}
 \qquad (x \in {\mathbb R}^m),
\nonumber\end{align}
from which Eq.\
$(\ref{eeformula})$
can be read off.
\phantom{.}\hfill$\qed$\par\noindent
\par\medskip\noindent
In the special case
$m = 2,$
working out the repeated binomial sums using Gould's combinatorial identities
(3.48) and (3.47) \cite{rGould} gives
\begin{align}\label{etwodim}
 e_{n_1, n_2} &= \pmatrix{2 \cr n_2 - n_1 + 1}
 = \cases{2 & \hbox{\rm if $n_1 = n_2$}, \cr 1 & \hbox{\rm if $|n_1 - n_2| = 1$}, \cr
 0 & \hbox{\rm otherwise}.}
\end{align}
This reflects the obvious fact that with only two types of objects, the
non-adjacency condition enforces an alternating arrangement, for which there
are two possibilities if the numbers of objects of both types are equal,
one possibility if they differ by one, and no possibility otherwise.
We emphasise that
$m = 2$
is a rather untypical case and that for
$m \ge 3$
much more complex arrangements are possible.
\par\medskip\noindent
{\it Remark.\/}
Eq.\
$(\ref{eeformula})$
can be given a direct combinatorial interpretation (and a somewhat more
convoluted proof) in the following manner.
We call any arrangement of the
$|n|$
objects ignoring the non-adjacency condition with
$t$
objects that are each followed by an object of the same type marked with a tick
an
{\it annotated arrangement with $t$ ticks.\/}
(Clearly there is no annotated version with
$t$
ticks of any arrangement which has fewer than
$t$
objects followed by an object of the same type.)
For each
$t \in {\mathbb N}_0,$
let
$A_t$
be the set of all annotated arrangements with
$t$
ticks. Its cardinality is
\begin{align}
 |A_t| &= \sum_{k \in {\mathbb N}_0^m, |k| = t} \pmatrix{n-1_m \cr k} \pmatrix{|n-k| \cr n-k}.
\label{enAt}\end{align}
Indeed,
given any element of
$A_t,$
we find
$k_j$
ticked objects of type
$j \in \{1, \dots, m\},$
so
$|k| = t.$
Considering the
$n_j$
objects of type
$j$
in the arrangement (ignoring the other types for the moment), the ticks can
occur in
$n_j - 1$
places, so there are
$\pmatrix{n_j-1 \cr k_j}$
possibilities.
Taking ticked objects together with their following object and single unticked
objects as groups, there will be
$n_j - k_j$
such groups. Among all types, the groups can be arranged in
$\pmatrix{|n-k| \cr n-k}$
ways, hence we obtain Eq.\
$(\ref{enAt}).$
Now to verify Eq.\
$(\ref{eeformula}),$
consider an arrangement of the objects.
Let
$l \in \{0, \dots, |n| - m\}$
be the number of objects in this arrangement followed by an object of the same
type. Annotated versions of this arrangement will appear in the sets
$A_0, \dots, A_l.$
In the set
$A_t,$
it will have
$t$
ticks which can be placed in
$l$
places, so there are
$\pmatrix{l \cr t}$
annotated versions of this (fixed) arrangement in this set.
We now count the total of its appearances (as different annotated arrangements) in
the sets
$A_0, \dots, A_l,$
counting its appearances in odd-indexed sets negative, those in even-indexed
sets positive. Thus in total we count this arrangement
\begin{align}
 \sum_{t=0}^l (-1)^t \pmatrix{l \cr t} &= (1-1)^l = \delta_{l, 0}
\nonumber\end{align}
times.
Hence the only arrangements counted in the alternating total are those with
$l = 0$
objects followed by an object of the same colour. This gives
\begin{align}
 e_n &= \sum_{t=0}^{|n|-m} (-1)^t |A_t|
\nonumber\end{align}
and hence, by Eq.\
$(\ref{enAt}),$
formula
$(\ref{eeformula}).$
\vskip 30pt

\section{The Number of Joint Ordered Factorisations}\label{s3}
\par\medskip\noindent
Given an $m$-tuple of integers
$a = (a_1, a_2, \dots, a_m) \in {\mathbb N}^m,$
with
$a_j \ge 2$
$(j \in \{1, \dots, m\}),$
we can use Theorem \ref{tcolours} to count the joint ordered factorisations of
$a$
where
$a_j$
is split into a prescribed number
$n_j$
of non-trivial factors. Indeed, we can think of taking
$n_j$
placeholders marked as type
$j$
(and otherwise indistinguishable), for
$j \in \{1, \dots, m\},$
and arranging all these placeholders according to the rules of the auxiliary
counting problem of Section 3, and then putting the factors for each
$a_j$
into the blocks of type
$j$
in their given order.
As there are
$e_n$
admissible arrangements of the placeholders and
$c_{n_j}(a_j)$
different non-trivial ordered factorisations of
$a_j,$
we obtain the number of joint ordered factorisations
from Eq.\
$(\ref{eeformula})$
as
\begin{align}
 e_{n} \prod_{j=1}^m c_{n_j}(a_j)
 &= \sum_{0_m \le k \le n-1_m} (-1)^{|k|} \pmatrix{n - 1_m \cr k} \pmatrix{|n-k| \cr n-k}
 \prod_{j=1}^m c_{n_j}(a_j).
\label{ecountjofr}\end{align}
The sum over all
$n \in {\mathbb N}^m$
(which is a finite sum since
$c_{n_j}(a_j) = 0$
if
$n_j$
exceeds
$\Omega(a_j),$
the number of prime factors of
$a_j$
counting multiplicities) then gives the total number of joint ordered
factorisations of
$a,$
which can be expressed as follows.
\begin{Theorem}\label{tcountjof}
Let
$m \in {\mathbb N}$
and
$a \in {\mathbb N}^m$
such that
$a_j \ge 2$
$(j \in \{1, \dots, m\}).$
Then the number of different joint ordered factorisations of
$a$
is
\begin{align}
 N_a &= \sum_{l \in {\mathbb N}^m} \pmatrix{|l| \cr l} \prod_{j=1}^m c_{l_j}^{(-l_j)}(a_j).
\label{eNa}\end{align}
\end{Theorem}
\par\medskip\noindent
{\it Proof.\/}
Summing the expression
$(\ref{ecountjofr})$
over
$n \in {\mathbb N}^m,$
changing the summation variable by setting
$n = k+l$
and interchanging the order of the sums, we find
\begin{align}
 N_a &= \sum_{k \in {\mathbb N}_0^m} \sum_{l \in {\mathbb N}^m} \pmatrix{k + l - 1_m \cr k} \pmatrix{|l| \cr l} \prod_{j=1}^m c_{k_j + l_j}(a_j)
\nonumber\\
 &= \sum_{l \in {\mathbb N}^m} \pmatrix{|l| \cr l} \prod_{j=1}^m \left(\sum_{k=0}^\infty \pmatrix{k+l_j-1 \cr k} c_{k+l_j}(a_j) \right).
\nonumber\end{align}
Now consider the functions which appear in the right-hand side product.
By the binomial formula for negative powers in the Dirichlet convolution
algebra of arithmetic functions,
\begin{align}
 \sum_{k=0}^\infty &(-1)^k \pmatrix{k + l_j-1 \cr k} c_{l_j+k}
 = (1-e)^{*l_j} * \sum_{k=0}^\infty \pmatrix{k + l_j - 1 \cr k} (e - 1)^{*k}
\nonumber\\
 &= (1-e)^{*l_j} * (e - (e - 1))^{*-l_j}
 = (1-e)^{*l_j} * \mu^{*l_j}
 = c_{l_j}^{(-l_j)},
\nonumber\end{align}
and hence the result.
\hfill$\qed$\par\noindent
\par\medskip\noindent
{\it Remarks.\/}
1.\ In view of the interpretation given to Eq.\
$(\ref{especial})$
in Section 2, Eq.\
$(\ref{eNa})$
can be read as
\begin{align}
 N_a &= (-1)^{\sum_{j=1}^m \Omega(a_j)} \sum_{n=m}^\infty (-1)^n \sum_{l \in {\mathbb N}^m, |l| = n} \pmatrix{|l| \cr l} \prod_{j=1}^m F_{l_j}(a_j),
\nonumber\end{align}
where we denote by
$F_{l_j}(a_j)$
the number of ordered factorisations of
$a_j$
into
$l_j$
non-trivial, square-free factors. Thus
$N_a$
can be construed as an alternating sum over
$n$
of the number of ways the integers
$a_1, a_2, \dots, a_m$
can be split into a total of
$n$
non-trivial, square-free factors, and these factors can be linearly arranged
without further constraints.

2.\ In the two-dimensional case $m=2$, Eq.\ $(\ref{etwodim})$ gives a simple explicit form for $e_{(n_1,n_2)}$,
and we obtain directly from the left-hand side of Eq.\ $(\ref{ecountjofr})$ that
\begin{align}
 N_{(a_1, a_2)} &= \sum_{n=1}^\infty 2\, c_n(a_1)\, c_n(a_2) + \sum_{n=1}^\infty c_n(a_1)\,c_{n+1}(a_2) + \sum_{n=1}^\infty c_{n+1}(a_1)\,c_n(a_2).
\nonumber\end{align}
In the symmetric case
$a_1 = a_2 = a$
considered in \cite{rOB,rHHLS}, this gives, via
$c_n + c_{n+1} = (1-e)^{*n}*(e + 1 - e) = c_n^{(1)}$,
the expression
$N_{(a, a)} = 2 \sum\limits_{n=1}^\infty c_n(a)\,c_n^{(1)}(a)$,
which involves different divisor functions compared to Eq.\ $(\ref{eNa})$ and reproduces Theorem 4 of
\cite{rHHLS}. Note that the count given in \cite{rHHLS} is $N_{(a, a)}/2$, using the permutation
symmetry of the two equally-sized component sets of the sum system.
More generally, when we have $m$ equally sized component sets in the sum system, then by the same permutation
symmetry,
$N_{(a, a, \dots, a)}/m!\in\mathbb{N}_0$.
Clearly, this property extends to more general $m$-tuples $a\in\mathbb{N}^m$ provided that all numbers
$a_j$
have the same factorisation structure, i.e.\ the multisets of exponents in the prime factorisation coincide.
The integer sequences
$(N_{(a, a)}/2!)_{a\in\mathbb{N}}$ and $(N_{(a, a, a)}/3!)_{a\in\mathbb{N}}$ are equal to sequences A0273013 and A0131514 in the OEIS (http://oeis.org), respectively, but it seems that no such OEIS record exists for $m\geq 4$.

We conclude with the observation that the number of $m$-part sum systems is at least $m!$, and this value is attained if and only if all parts have prime cardinality.
\begin{Theorem}
Let
$m \in {\mathbb N}$
and
$a \in {\mathbb N}^m$
such that
$a_j \ge 2$
$(j \in \{1, \dots, m\}).$
Then
$N_a\ge m!$
and equality holds if and only if all
$a_j$
are prime numbers.
\end{Theorem}
\noindent
{\it Proof.\/}
Starting from the left-hand side of Eq.\ $(\ref{ecountjofr})$ and considering that
$e_{1_m} = m!$ by Eq.\ $(\ref{eeformula})$ and $c_1(a_j) = 1$, we find
\begin{align}
 N_a &= m! + \sum_{n\in\mathbb{N}\setminus\{1_m\}} e_n \prod_{j=1}^m c_{n_j}(a_j) \ge m!.
\label{eprime}
\end{align}
Suppose one of the dimensions, w.l.o.g.\ $a_1$, is not a prime, and consider $n = (2, 1, \dots, 1)$.
Then
$\prod\limits_{j=1}^m c_{n_j}(a_j) = c_2(a_1) \ge 1$
and
$e_{(2, 1, \dots, 1)} = \frac{(m+1)!}{2!}$ by Eq.\ $(\ref{eeformula})$,
so the sum in $(\ref{eprime})$ is strictly greater than 0.
\hfill$\qed$
\section*{Appendix}
\par\medskip\noindent
For the reader's convenience, here is a summary of the standard multi-index
notation used in this paper.
For $m$-tuples of integers, we apply the
usual componentwise addition, subtraction and scalar multiplication as well as the
size function and the multi-factorial,
\begin{align}
 |n| &= \sum_{j=1}^m |n_j|, \qquad
 n! = \prod_{j=1}^m n_j! \qquad (n \in {\mathbb N}_0^m),
\nonumber\end{align}
respectively, and the partial ordering
\begin{align}
 n \le \tilde n &\Longleftrightarrow n_j \le \tilde n_j \qquad (j \in \{1, \dots, m\}).
\nonumber\end{align}
We define the special $m$-tuples
$0_m = (0, \dots, 0), 1_m = (1, \dots, 1) \in {\mathbb N}_0^m.$
In addition to the usual binomial coefficients, we use the multi-binomial
coefficients
\begin{align}
 \pmatrix{n \cr \tilde n} &= \prod_{j=1}^m \pmatrix{n_j \cr \tilde n_j} = \frac{n!}{\tilde n!\,(n - \tilde n)!} \qquad
 (n, \tilde n \in {\mathbb N}_0^m, \tilde n \le n)
\nonumber\end{align}
and the multinomial coefficients
\begin{align}
 \pmatrix{|n| \cr n} &= \frac{|n|!}{\prod_{j=1}^m n_j!} = \frac{|n|!}{n!} \qquad (n \in {\mathbb N}_0^m).
\nonumber\end{align}
Note that using the same bracket notation for these different quantities
does not create confusion, since the type (dimensionality) of the arguments
determines which coefficient is meant.
Finally,
\begin{align}
 x^n &= \prod_{j=1}^m x_j^{n_j} \qquad (n \in {\mathbb N}_0^m, x \in {\mathbb R}^m).
\nonumber\end{align}


\begin{thebibliography}{99}
\bibitem{rGEA}
G.E.~Andrews.
{\it Number Theory.\/}
Dover, New York 1994
\bibitem{rGould}
H.W.~Gould.
{\it Combinatorial identities.\/}
Morgantown 1972
\bibitem{rHHLS}
S.L.~Hill, M.N.~Huxley, M.C.~Lettington, K.M.~Schmidt.
Some properties and applications of non-trivial divisor functions.
{\it Ramanujan J. \/}
(2019) https://doi.org/10.1007/s11139-018-0093-9
\bibitem{rHLS}
M.N.~Huxley, M.C.~Lettington, K.M.~Schmidt.
On the structure of additive systems of integers.
{\it Periodica Mathematica Hungarica\/}
{\bf 78}
(2019) 178--199
\bibitem{rIvic}
A.~Ivi\'c.
On the mean square of the divisor function in short intervals.
{\it J.\ Th\'eor.\ Nombres Bordeaux\/}
{\bf 21}
(2009) 251--261
\bibitem{rKRRR}
J.P.~Keating, B.~Rodgers, E.~Roditty-Gershon and Z.~Rudnick.
Sums of divisor functions in $F_q[t]$  and matrix integrals.
{\it Math.\ Z.\/}
{\bf 288}
(2018) 167--198
\bibitem{rOB}
K.~Ollerenshaw and D.~Br\'ee.
{\it Most-perfect pandiagonal magic squares.\/}
IMA 1998
\bibitem{rPopo}
C.P.~Popovici.
O generalizare a func\c tiei lui M\"oebius.
{\it Acad.\ R.\ P.\ Rom\^\i ne Stud.\ Cerc.\ Mat.\ \/}
{\bf 14}
(1963) 493--499
\bibitem{rRama}
S.~Ramanujan.
On the number of divisors of a number.
{\it J.\ Indian Math.\ Soc.\/}
{\bf 7}
(1915) 131-133
\bibitem{rSchSp}
W.~Schwarz and J.~Spilker.
{\it Arithmetical Functions.\/}
{\it LMS Lecture Note Series\/}
{\bf 184.}
Cambridge University Press, Cambridge 1994
\bibitem{rStan}
R.P.~Stanley.
{\it Enumerative Combinatorics, Vol.\ 1.\/}
Cambridge Univ.\ Press, Cambridge 1997
\end{thebibliography}
\end{document}